\documentclass[12pt]{amsart}
\usepackage{amsmath}
\usepackage{amssymb, amstext, amsthm, amsfonts,euscript,amscd,mathrsfs}%accents}
\usepackage[colorlinks=false]{hyperref}
\usepackage{graphicx}
\usepackage[final]{pdfpages}
\usepackage[toc,page]{appendix}
\usepackage{latexsym}
\usepackage[disable]{todonotes}

\DeclareMathOperator{\Real}{\mathbb{R}}

\newcommand{\N}{\mathbb{N}}

\newcommand{\jap}[1]{\!\left<#1\right>}
  \newcommand{\jxi}{\jap{\xi}}

\newcommand{\dx}{\partial_x}
\newcommand{\dxp}[1]{\partial_{x_{#1}}}

\newcommand{\dy}{\partial_y}
\newcommand{\dyp}[1]{\partial_{y_{#1}}}

\providecommand{\norm}[1]{\lVert#1\rVert}
\providecommand{\abs}[1]{\lvert#1\rvert}

\newcommand{\conv}{\ast}

\DeclareMathOperator{\dom}{dom}

\newcommand{\tld}[1]{\tilde{#1}}

\newcommand{\R}{\Real}

\newcommand{\eps}{\varepsilon}

\newcommand{\Hormander}{{H\"{o}rmander}}

  \DeclareMathOperator{\supp}{supp}
   
\newcounter{theorem_counter}
\newtheorem*{thm*}{Theorem}
\newtheorem{thm}{Theorem}
\newtheorem{lem}[thm]{Lemma}
\newtheorem{prop}[thm]{Proposition}
\newtheorem{rem}[thm]{Remark}
\newtheorem{coro}[thm]{Corollary}

\newtheorem{defn}[thm]{Definition}

\newtheorem{conj}[thm]{Conjecture}

\title[Logarithmic estimate for hypoellipticity]{Necessity of a logarithmic estimate for hypoellipticity of some degenerately elliptic operators}
\author{Timur Akhunov}\address{Department of Mathematics and Computer Science\\
Wabash College\\
301 W Wabash Ave, Crawfordsville, IN 47933, USA}
\author{Lyudmila Korobenko}
\address{Mathematics Department\\
Reed College\\
3203 Southeast Woodstock Boulevard\\
Portland, OR 97202-8199, USA}

\thanks{The authors thank Cristian Rios, who introduced us to the problem of hypoellipticity and worked with us on \cite{AkhKorRios} that motivated results of this paper. We also thank the anonymous referee for numerous insightful suggestions}
\date{\today}
\begin{document}
\begin{abstract}
  This paper extends a class of degenerate elliptic operators for which hypoellipticity requires more than a logarithmic gain of derivatives of a solution in every direction. Work of Hoshiro and Morimoto in late 80s characterized a necessity of a super-logarithmic gain of derivatives for hypoellipticity of a sum of a degenerate operator and some non-degenerate operators like Laplacian. The operators we consider are similar, but more general. We examine operators of the form $L_1(x)+g(x)L_2(y)$, where $L_1(x)$ is one-dimensional and $g(x)$ may itself vanish. The argument of the paper is based on spectral projections, analysis of a spectral differential equation, and interpolation between standard and operator-adapted derivatives. Unlike prior results in the literature, our methods do not require explicit analytic construction in the non-degenerate direction. In fact, our result allows non-analytic and even non-smooth coefficients for the non-degenerate part.
\end{abstract}
\maketitle
%\section{Old Introduction - Rewrite completely}
\section{Introduction}
Consider a degenerate elliptic operator %\footnote{defined in \eqref{L2} and \eqref{elliptic} for $L_2$ and similarly for $L_1$}
\begin{align}\label{Lgeneral}
  L(x,y)=L_1(x)+g(x)L_2(y)
\end{align}
where
\begin{align*}
  L_1=-\sum_{j,k=1}^n a_{jk}(x)\dxp{j}\dxp{k}  +\sum_{j=1}^n a_j(x) \dxp{j} + a_0(x)
\end{align*}
with $x\in \R^n$, $a_{jk}$, $a_j$ smooth real valued coefficients with $a_{jk}$ a non-negative matrix. Similarly, let $L_2$ be denoted as
\begin{align}
  \label{L2}
  L_2 = -\sum_{j,k=1}^m b_{jk}(y)\dyp{j}\dyp{k} +\sum_{j=1}^m b_j(y) \dyp{j} + b_0(y)
\end{align}
where $y\in \R^m$,  $b_{jk}$, $b_j$ are smooth real functions and $b_{jk}$ is a non-negative matrix:
\begin{align}\label{elliptic}
  \sum_{j,k=1}^m b_{j,k}(y) \eta_j \eta_k \ge 0
\end{align}
In this paper we investigate conditions under which local hypoellipticity of $L$, or smoothness of solutions of $Lu=f$ whenever $f$ is locally smooth (see Definition \ref{def-hypo}), requires the following superlogarithmic estimate introduced by Morimoto \cite{Mor87}:
 \begin{defn}[Superlogarithmic estimate]\label{def-superlog}
   We say that an operator $L$ satisfies a superlogarithmic estimate near $y_0\in R^m$ if for any small enough compact set $K\Subset R^m$ containing $y_0$ and for all $\eps>0$, there exists a family of constants $C_{\eps,K}$, such that the following estimate holds.
\begin{equation}\label{superlog}
||\log\jap{\xi}\hat u(\xi)||^2 \leq \varepsilon\,\mathrm{Re}(L u,u) + C_{\varepsilon,K} ||u||^2,\quad u\in C_{0}^{\infty}(K).
\end{equation}
where
 \begin{align}\label{jap-bracket}
   \jap{\xi}:=\sqrt{e^2+|\xi|^2}
 \end{align}
 \end{defn}
This estimate is satisfied by certain operators $L_2$ with degeneracy in \eqref{elliptic} so that the vanishing is faster than any polynomial in $y$ in some directions. For example,
\begin{align}\label{Kusuoka}
  L_2=\partial_{y_1}^2 + e^{-\frac{1}{|y_1|^{1-\kappa}}}\partial_{y_2}^2
\end{align}
considered by \cite{Kusuoka-Strook84} satisfies \eqref{superlog} for $\kappa>0$, see \cite{Christ01} Lemma 5.2. The {\Hormander} bracket condition for $L_2$ corresponds to a polynomial degeneracy for $L_2$ and leads to subelliptic estimate, defined as
 \begin{align}\label{subelliptic}
||\jap{\xi}^\delta \hat u(\xi)||^2 \leq C(\delta,K)\left(\mathrm{Re}(L_2 u,u) + ||u||^2\right),
   \text{ for }\quad u\in C_{0}^{\infty}(K).
 \end{align}
 for some $\delta>0$ and $\jap{\xi}$ defined in \eqref{jap-bracket} as before. Note that \eqref{subelliptic} implies \eqref{superlog}. Superlogarithmic estimate \eqref{superlog} limits the degree of vanishing that an operator $L_2$ can exhibit. If a stronger vanishing is considered, e.g. \eqref{Kusuoka} for $\kappa\le0$, then \eqref{superlog} fails and a weaker gain function should replace the logarithm in \eqref{superlog}. See \cite{Christ01} for a further discussion of degeneracy and gain function.\\

  This is a companion paper to \cite{AkhKorRios}, where we proved hypoellipticity of $L$ from \eqref{Lgeneral} under the assumption that operators $L_1$ and $L_2$ satisfy Definition \ref{def-superlog}, see \cite[Theorem 1]{AkhKorRios} and \cite{Mor87}. It was also demonstrated that, at least for particular classes of operators, the superlogarithmic estimate \eqref{superlog} for $L_1$ is necessary, see \cite[Example 1]{AkhKorRios}. Here we address in greater generality the question of necessity of these superlogarithmic estimates. \\

  We first set up relevant notation. In what follows, for any subsets $E$, $F$ in $\R^n$ we write $E\Subset F$ to indicate that the closure of $E$, $\overline{E}$, is compact and $\overline{E}\subset F$.\\

  For completeness, we first define local smooth hypoellipticity here.
       \begin{defn}\label{def-hypo}
  The operator $L(x,y)$ is {\bf locally ($C^\infty$) hypoelliptic at $(x_0,y_0)$} if there exists a neighborhood $\Omega$ of $(x_0,y_0)$, such that for every $f(x,y)$ smooth in $\Omega$ and a distributional solution $u$ satisfying $Lu=f$ in the weak sense in $\Omega$, then $u\in C^\infty(\Omega')$ for all neighbohoods $\Omega'$ of $(x_0,y_0)$ such that $\Omega'\Subset \Omega$.
\end{defn}
In particular, the region of smoothness of a weak solution depends on the operator $L$ and the region $\Omega$, but not on the specific function $f$, so long as $f$ is smooth. Our method allows us to work with a limited regularity form of hypoellipticity as well. We model this form of hypoellipticity on the Definition \ref{def-hypo}.
\begin{defn}\label{hypo-weak}
  We say that $L$ is {\bf (locally) $H^s$-hypoelliptic at $(x_0,y_0)$}, if there exists a neighborhood $\Omega$ of $(x_0,y_0)$, such that for every $f\in H^s(\Omega)$ and $u\in L^2(\Omega)$ satisfying $Lu=f$ in the weak sense in $\Omega$, then $u\in H^s(\Omega')$ for all neighbohoods $\Omega'$ of $(x_0,y_0)$ such that $\Omega'\Subset \Omega$.
\end{defn}

  Morimoto in \cite{Mor87}, motivated by the results in \cite{Kusuoka-Strook84}, provided a non-probabilistic proof that in order for a symmetric operator $L(x,y)=-\dx^2 + L_2(y)$, $x\in \R,\  y\in\R^m$ to be hypoelliptic, $L_2$ has to satisfy Definition \ref{def-superlog}. Independently at the same time, Hoshiro \cite{Hos87} obtained a similar result as a corollary of an a priori estimate for a general self-adjoint partial differential operator.
Our main theorem is a generalization of the Morimoto's and Hoshiro's necessity estimates to the operator of the form $L=L_1(x)+g(x)L_2(y)$, where the function $g$ is allowed to vanish, allowing the operator $L$ itself to be more degenerate than $L_2$ and hence not satisfy \eqref{superlog} at some $(x_0,y_0)$.
 \begin{thm}\label{thm:main}
 Consider a possibly degenerate elliptic operator
	\begin{align}\label{tldL}
	L = -a_2(x)\dx^2 + a_1(x)\dx + a_0(x) + g(x)L_2(y,\dy),
	\end{align}
	where $L_2$, defined in \eqref{L2}, is a self-adjoint operator, i.e.
   \begin{align}\label{symmetry}
  (L_2 u, v) = (u, L_2 v) \text{, for }  u\in C^\infty_0
\end{align}
Assume further $a_2(x_0)> 0$, $g(x)\ge 0$ and the coefficients $g$ and $a_i$ for $i=0,1,2$ are smooth.\\

Suppose that $L$ from \eqref{tldL} is {\bf $H^s$-hypoelliptic} near $(x_0,y_0)$ for {\bf $s>\frac{1}{2}$}. Then $L_2$ must satisfy \eqref{superlog} near $y_0$.
\end{thm}
\begin{rem}
  Proof of Theorem \ref{thm:main} for $\frac{1}{2}<s\le 1$ is valid with $a_2(x)\in C^1_x$ and other coefficients merely continuous.
\end{rem}
%We remark that the proof of this Corollary follows the same argument as that of the Theorem \ref{thm:main}. The need for continuity of the coefficients is particularly evident in Lemma \ref{lem:ODE}.\\%, such regularity of a solution does not require more than continuity of , unless more is needed for $H^s$ hypoellipticity, which is a black box for our result. Thus an analogue of Theorem \ref{thm:main} is possible even with limited regularity of coefficients.

\begin{coro}\label{coro:smooth}
	Suppose $L$ from \eqref{tldL} satisfies all of the assumptions of Theorem \ref{thm:main}, with the exception of $H^s$  hypoellipticity. We replace that assumption with the $C^\infty$ hypoellipticity of $L$ at $(x_0,y_0)$. Then  $L_2$ satisfies \eqref{superlog} near $y_0$. Note that to establish $C^\infty$ hypoellipticity, more smoothness of the coefficients than continuity in Theorem \ref{thm:main} may be needed.
\end{coro}

\begin{rem}\label{rem:g-vanishing}
	Note, that the focus in Theorem \ref{thm:main} is on the implications of hypoellipticity. In particular, if $g(x)\equiv 0$, the operator $L$ fails to be hypoelliptic and the Theorem \ref{thm:main} does not apply.\\% even if $L_2$ was elliptic.\\

When $g(x)$ has an isolated zero at $x_0$, our earlier work with Cristian Rios in \cite{AkhKorRios} proved $C^\infty$  hypoellipticity of $L$ if it satisfies \eqref{superlog} extending the non-degenerate case of $g(x)$ studied by Morimoto in \cite{Mor87}. Thus, in the context of an isolated zero (or non-vanishing) $g$, Theorem \ref{thm:main} is sharp for $C^\infty$ hypoellipticity of \eqref{tldL}.\\

Finally, \cite{Morimoto-Morioka97-Fedii} proved hypoellipticity of $L$ from \eqref{tldL}, when $L_2=\dy^2$ assuming only that the stochastic average
\begin{align*}
  g_I:=\int_I g(x) dx>0 \text{, for any interval }I\subset \R,
\end{align*}
including the case of non-isolated zeros. Morimoto and Morioka further connected quantitative estimates for $g_I$ with the logarithmic estimate \eqref{superlog} for some operators in \cite{Morimoto-Morioka97}. Thus Theorem \ref{thm:main} is likely non-trivial in some cases of $g(x)$ with non-isolated zeros.
%  $g(x)$ has more complicated vanishing structure than an isolated zero at $x_0$ and is not identically $0$, proving hypoellipticity
%
%  no assumption is made on $g(x)$ in the Theorem \ref{thm:main}, except the hypothesis that $L$ is hypoelliptic. For, example if $g(x)\equiv 0$ and $y\in R^n$ for $n\ge 1$, $L$ will fail to be hypoelliptic regardless of whether \eqref{superlog} holds for $L_2$. However, authors' earlier work \cite{AkhKorRios} shows that if $x_0$ is an isolated zero of $g(x)$ and $g(x)\ge 0$, $L$ must be hypoelliptic if it satisfies \eqref{superlog}. In that context, \eqref{superlog} is sharp for hypoellipticity of \eqref{tldL}.
\end{rem}
The strength of Theorem \ref{thm:main} lies both in the fact that we replace $-\partial^{2}_{x}$ by a general one-dimensional elliptic operator, and, more importantly, allow for broader conditions on the degeneracy of the non-negative weight function $g(x)$ in \eqref{tldL}. This is in line with the results by Fedi\u{\i} \cite{Fedii71}, who considered $L=\partial^{2}_{x}+g(x)\partial^{2}_{y}$, and Kohn \cite{Kohn98} who generalized Fedi\u{\i}'s result to $L=L_1(x)+g(x)L_2(y)$ where $L_1$ and $L_2$ are subelliptic in their variables. Note that the subelliptic estimate \eqref{subelliptic} is stronger that (\ref{superlog}), so Theorem \ref{thm:main} together with \cite[Theorem 1]{AkhKorRios} characterizes hypoellipticity in the case of $g(x)$ with isolated vanishing. Furthermore, Theorem \ref{thm:main} extends a necessary condition for hypoellipticity to a wide class of operators, many of them not previously considered in the literature.\\

Theorem \ref{thm:main} can be generalized in several ways. First, one can ask if the result still holds for a more general operator $L_1(x)$. We believe that the ultimate goal would be the following conjecture.
\begin{conj}\label{conj}
  Suppose that the operator $L$ from \eqref{Lgeneral} is symmetric, i.e. satisfies \eqref{symmetry} and $g(x)\ge 0$. The operator $L$ is locally $C^\infty$-hypoelliptic, if and only if both $L_1$ and $L_2$ satisfy \eqref{superlog}.
\end{conj}
As discussed in Remark \ref{rem:g-vanishing}, results in \cite{AkhKorRios} imply the sufficiency part of this conjecture in the case of $g(x)$ that has at most isolated zeros. Theorem \ref{thm:main} confirms the necessity direction in the case that $L_1$ is one dimensional and elliptic. \\%

 The second direction to generalize Theorem \ref{thm:main} is to allow for a more general dependence on the variables, ultimately considering operators of the form $L_1(x,y)+L_2(x,y)$ with either $L_1$ or $L_2$ being degenerate elliptic. It is a nontrivial question, what happens to hypoellipticity when two (or more) operators are added. For example, Fedi\u{\i}'s operator $\partial^{2}_{x}+g(x)\partial^{2}_{y}$ is hypoelliptic whenever $g$ is positive away from zero \cite{Fedii71}, while a very similar three-dimensional operator $\partial^{2}_{x}+\partial^{2}_{y}+g(x)\partial^{2}_{z}$ considered by Kusuoka and Strook \cite{Kusuoka-Strook84} is hypoelliptic if and only if $g$ vanishes slower than $e^{-1/x^2}$. See also \cite{Christ01}, \cite{Mor87}, \cite{Morioka92}, and references within for more discussion of these ``intertwined cases''. It is possible that methods used in our paper may also extend to some complex valued elliptic operators satisfying {\Hormander} bracket condition \cite{Kohn05}, similar to the note of Christ \cite{Christ05}.\\

 To summarize, Theorem \ref{thm:main} together with \cite[Theorem 1]{AkhKorRios},  contribute to our understanding of how regularity properties of operators in lower dimensions play into the regularity of their sum in higher dimensions. On the technical level, prior works of \cite{Hos87}, \cite{Morimoto-nonhypo} and \cite{Mor87} used an explicit construction of an analytic function of $x$ motivated by techniques of analytic hypoellipticity \cite{Metivier78} or complex analysis. In contrast, the argument of this paper is grounded in analysis of ODEs, which allows for non-analytic and even nonsmooth coefficients.
Additionally, we have adapted Littlewood-Paley projections to the degenerate operator to interpolate between standard derivatives and functional calculus, what we hope improves exposition.\\

This paper is organized as follows. The proof of Theorem \ref{thm:main} is split into the next three sections. Section \ref{sec:ODE:soln} uses an initial value problem for an ODE and spectral projections to construct a special solution to $\tilde{L}u=0$. In Section \ref{sec:projections} we use the Closed Graph Theorem and a version of Sobolev Embedding to obtain qualitative estimates for the solution of a hypoelliptic operator. In section \ref{projection:gain} we apply the results of the previous two section to spectral projections of any $L^2$ function. Section \ref{sec:interpolation} concludes the proof by using interpolation to combine the estimates for the spectral projections into (\ref{superlog}). Finally, the Appendix \ref{sec:spectral}, contains some technical results about spectral projections used throughout the paper.

%\hide{ i.e. Results of Hoshiro and Morimoto survive lower order perturbation}

    \section{Spectral solutions via an ODE}\label{sec:ODE:soln}
        We begin the proof of Theorem \ref{thm:main} with a simple change of operator. It is convenient to reduce \eqref{tldL} to a case of a constant coefficient: $a_2\equiv 1$ first.

\begin{lem}\label{rem:a2is1}
Consider the operator $\tld P = \frac{1}{a_2(x)} L$ in a neighborhood small enough, where $a_2(x)$ does not vanish. I.e.
  \begin{align*}
    \tld P = -\dx^2 + \frac{a_1}{a_2}(x)\dx + \frac{a_0}{a_2}(x) + \frac{g}{a_2}(x)L_2(y,\dy).
  \end{align*}
Then $\tld P$ is hypoelliptic if and only if $L$ from \eqref{tldL} is hypoelliptic. By relabeling $a_1$, $a_0$ and $g(x)$ with analogues in $\tld P$, we can assume $a_2(x)=1$ in $L$ without loss of generality.
  %Dividing $L$ by $a_2(x)\neq 0$ preserves hypoellipticity.
\end{lem}
\begin{proof}
  Suppose $L$ is hypoelliptic. Note that $\tld P$ has the same structure as $L$, except for a coefficient of $\dx^2$.\\

      If $F=\tld P u$ is smooth near $(x_0,y_0)$, so is $f(x,y)=a_2(x)F(x,y)=L u$. By hypoellipticity of $L$, $u$ must be smooth near $(x_0,y_0)$. Thus $\tld P$ is hypoelliptic. A similar argument works for the converse.
\end{proof}
%\comment{9/2022 To be complementary, if $a_2(x_0)=0$ a tweak of a $\delta$ function may create a solution ruining hypoellipticity. We choose not to do this}
The strategy of the proof of Theorem \ref{thm:main} is to convert the problem $L u= 0$ into an ODE using spectral projections, formally discussed in appendix \ref{sec:spectral}. Informally, the operator $L_2$ is replaced by a spectral parameter $\lambda$. The following ODE becomes relevant to our analysis in light of Lemma \ref{rem:a2is1}.% follow ODE (This section is concerned
\begin{align}\label{lambda-ode}
        \begin{cases}
          -\dx^2 v(x,\lambda) + a_1(x)\dx v(x,\lambda)+ [a_0(x) + g(x)\lambda]v(x,\lambda) =0 \\
          v(x_0,\lambda)=1,\, \dot v(x_0,\lambda)=0
        \end{cases}
       \end{align}
       We analyze this ODE first before returning to the original operator $L$ from \eqref{tldL} in Lemma \ref{lem:tld-spectral}.
  %  The spectral projection in the previous section converts the operator $L$ from \eqref{tldL} to an ODE \eqref{lambda-ode}. We use the following wellposedness result for a system of ODEs:
 \begin{prop}[Proposition 1.2.4, \cite{HorLec}]\label{ODE:prop}
   Let $x_0\in \R$, $I$ a connected neighborhood of $x_0$ and let $\Phi:I\to R^{m\times m}$ be a continuous function with values in $m\times m$ matrices. Then the initial value problem\begin{align}\label{ODE}
     \begin{cases}
     y'(x) = \Phi(x) y & \,\, \text{ for }|x-x_0|\le r_{x_0} \\ y(x_0)= y_0 &
   \end{cases}
   \end{align}
   has a unique solution for all $y_0\in \R^m$. Moreover, if $\norm{\Phi}_{L^\infty_I} \le M$ for some $M>1$, then $|y(x)|\le |y_0|e^{M|x-x_0|}$ for all $x\in I$.
 \end{prop}
 We need estimates for derivatives of the solutions to \eqref{ODE}, primarily for $y'(x)$.
 \begin{lem}\label{ODE:derivatives}
   If in addition $\norm{\Phi^{(j)}(x)} \le M$ for $0\le j\le k-1$, $k\ge 1$, then
   \begin{align*}
     |\dx^k y(x)|\le C(k)(1+M)^k e^{M|x-x_0|}|y_0|
   \end{align*}
   Moreover, estimates for the derivatives of $\Phi$ are not needed for $k=1$.
 \end{lem}
 \begin{proof}
   Note, that
 \begin{align*}
   |y'(x)|\le \abs{\Phi(x)}\abs{y(x)}\le M e^{M|x-x_0|}|y_0|
 \end{align*}
 Differentiating \eqref{ODE} we get
 \begin{align*}
   y''= \Phi y' + \Phi' y
 \end{align*}
 Hence
 \begin{align*}
   \abs{y''(x)}\le (M e^{M|x-x_0|}+ M \cdot M e^{M|x-x_0|})|y_0|
 \end{align*}
Proceeding inductively, we get the estimates for higher derivatives of $y$ up to $k$.
 \end{proof}
 We now use these results to solve \eqref{lambda-ode}.
 \begin{lem}\label{lem:ODE}
    Let $x_0$ be given. Then \eqref{lambda-ode} has a unique solution $v(x,\lambda)$ of \eqref{lambda-ode} on $[x_0-1,x_0+1]$. Furthermore, the following estimate holds:
    \begin{align*}
      |\dx^k v(x,\lambda)|\le C({k,x_0}) \jap{\lambda}^{\max\{0,\frac{k-1}{2}\}} \exp(C_{x_0}\jap{\lambda}^{\frac{1}{2}} |x-x_0|)
    \end{align*} % and depends continuously on $\lambda$.
    where $C_{x_0}$ depends on $x_0$ through the coefficients of \eqref{lambda-ode}. Moreover, the constant for $k=0$ and $k=1$, $C({0,x_0})=1$. For $k=0,1,2$ no derivatives of the coefficients $a_0$, $a_1$ and $g$ are needed. Note, that we use \eqref{jap-bracket} for $\jap{\lambda}$.
  \end{lem}
\begin{proof}
  To apply Proposition \ref{ODE:prop} to \eqref{lambda-ode}, define $a(x)= [a_0(x) + g(x)\lambda]$ and $y(x)=\begin{pmatrix}v \\ v'
 \end{pmatrix}$. Then $y(x_0)=\begin{pmatrix}1 \\ 0
 \end{pmatrix}$ and $y'(x) = \Phi(x) y(x)$ with
 \begin{align}\label{ODE-system}
   \Phi(x) = \begin{pmatrix}0 & 1\\a(x) & a_1(x)
 \end{pmatrix}
 \end{align}
 Looking at the eigenvalues of $\Phi(x)$, we obtain $$\norm{\Phi(x)} \le   C(|a_1(x)|+ \sqrt{|a(x)|}),$$ where $C$ is a universal constant.
 Therefore, from the definition of $a(x)$, we estimate
% For $$$$ small enough, $g(x)\approx_{x_0} 1$, so
 \begin{align}\label{Phi-norm}
   \norm{\Phi(\cdot)}_{L^\infty_{loc}} \le C(\sqrt{\norm{g}_{L^\infty_{loc}}\lambda+\norm{a_0}_{L^\infty_{loc}}} + \norm{a_1}_{L^\infty_{loc}})\le C(x_0)\jap{\lambda}^{\frac{1}{2}}
 \end{align}
 Proposition \ref{ODE:prop} completes the proof for $k=0$ and $k=1$.\\

 For higher derivatives of $v$, we observe that differentiating \eqref{ODE-system}, will give an estimate similar to \eqref{Phi-norm} for $\norm{\dx^j\Phi(\cdot)}_{L^\infty_{loc}}$ . That is, $\dx a(x)$ and $\dx a_1(x)$ have similar to dependence on $\lambda$, i.e. $$\norm{\dx^j\Phi(\cdot)}_{L^\infty_{loc}}\le C(j,x_0)\jap{\lambda}^{\frac{1}{2}}.$$
 As control of $v^{(k)}(x)$ requires estimates of $y^{(k-1)}(x)$ for $k\ge 1$, Lemma \ref{ODE:derivatives} completes the proof.
%We conclude with
\end{proof}

     % Sketch of the proof of Lemma is below the outline
 We are now ready to construct solutions to \eqref{tldL} via Proposition \ref{spectral} and other results from the Appendix. We define an operator $v(x,B)$ on a subset of $L^2$ via $v(x,\lambda)$. The following Lemma summarizes the connection between the spectral construction and \eqref{tldL}. We will later extend the domain of the $v$ operator, by pre-applying Littlewood-Paley projections.
 \begin{lem}\label{lem:tld-spectral}
Let $u(y)$ be such that, $e^{\sqrt{B}\eps}u\in L^2$ for some $\eps>0$. Define $w$ by
 \begin{align}\label{spectral-soln}
   w(x,y) = v(x,B) u(y) = \int_1^\infty v(x,\lambda) dE_\lambda u(y).
 \end{align}
 where $v(x,\lambda)$ is the solution of \eqref{lambda-ode}.
Let $0<\eps'<\eps$. Then for $|x-x_0|\le \frac{\eps'}{2 C_{x_0}}$ for $C_{x_0}$ from Lemma \ref{lem:ODE}, $w$ is well-defined. Likewise, $L w(x,y) \in L^2$ for $L$ from \eqref{tldL}. Moreover, using Lemma \ref{rem:a2is1} to replace $a_2$ with $1$, we establish that $$Lw =0.$$
 \end{lem}
 \begin{proof}
   By Proposition \ref{spectral}
 $$
 B w = \int_1^\infty \lambda v(x,\lambda) dE_\lambda u(y)
 $$
 Hence by Lemma \ref{lem:ODE} for $|x-x_0|\le \frac{\eps'}{2C_{x_0}}$
 \begin{align*}
   \norm{ B w}_{L^2}^2& \le \int_1^\infty  \lambda^2(\exp(C_{x_0}\jap{\lambda}^{\frac{1}{2}} |x-x_0|))^2  d \norm{E_\lambda u(y)}^2 \\
   & \le C\int_1^\infty \lambda^2 e^{\eps'\sqrt{ \lambda}} d \norm{E_\lambda u(y)}^2\le C\int_1^\infty e^{\eps\sqrt{ \lambda}}d \norm{E_\lambda u(y)}<\infty
 \end{align*}
 for $\eps'<\eps$ and $x$ within the interval specified above.\\
 We now justify that $\dx w\in L^2$ and can be expressed as
 \begin{align}\label{ODE:dx-w}
   \dx w = \int_1^\infty \dx v(x,\lambda) dE_\lambda u(y).
 \end{align}
 To do that we use difference quotients to interchange derivatives and integration in $\lambda$. Indeed,
\begin{align*}
  D_h w(x):=\frac{w(x+h,y)-w(x,y)}{h}= \int_1^\infty \frac{1}{h}[v(x+h,\lambda) -v(x,\lambda)]dE_\lambda u(y)
\end{align*}
Hence
\begin{align*}
  D_hw(x)=  \int_1^\infty \int_0^1 \dx v(x+sh,\lambda) ds dE_\lambda u(y)
\end{align*}
Since $|x-x_0|\le \frac{\eps'}{2 C_{x_0}}$ and for $|h|<\frac{\eps-\eps'}{2sC_{x_0}}$ we use Lemma \ref{lem:ODE} to obtain
\begin{align*}
  |\dx v(x+sh,\lambda)|\le \jap{\lambda}^{\frac{1}{2}} e^{C_{x_0}\eps\sqrt{ \lambda}|x-x_0+sh|} \le e^{\eps\sqrt{ \lambda}}
\end{align*}
Therefore, $D_hw\in L^2$ for $h$ sufficiently small. A similar argument shows that for any $h_n\to 0$, $D_{h_n} w$ is Cauchy and thus has a limit, that we call $\dx v$. By Lebesgue Dominated Convergence $\int_0^1 \dx v(x+sh,\lambda) ds \to v(x,\lambda)$ as $h\to 0$. This proves \eqref{ODE:dx-w}.\\

Analysis of $\dx^2 w\in L^2$ and an analogue of \eqref{ODE:dx-w} is similar.  Combining the terms implies $L w\in L^2$.\\

We now demonstrate that under the hypothesis of the Lemma $Lw =0$. Using spectral projections, e.g. Proposition \ref{spectral}
\begin{align*}
  L w(x,y) = -\dx^2 w + a_1(x)\dx w + \int_1^\infty [a_0(x) + g(x)\lambda]v(x,\lambda) dE_\lambda u(y)
\end{align*}
We have shown above that we can pass $x$ derivatives into the projection for $w$. Hence
$$L w(x,y) = \int_1^\infty [-\dx^2 v(x,\lambda) + a_1(x)\dx v(x,\lambda)+ [a_0(x) + g(x)\lambda]v(x,\lambda) dE_\lambda u(y)$$

  As $v(x,\lambda)$ solves \eqref{lambda-ode}, $L w=0$.
%$$L w(x,y) = \int_1^\infty [-\dx^2 v(x,\lambda) + a_1(x)\dx v(x,\lambda)+ [a_0(x) + g(x)\lambda]v(x,\lambda) dE_\lambda u(y)$=0$
%
%  As $L w =0$, from the hypoellipticity $w\in C^\infty_{loc}$. In particular, as from \eqref{lambda-ode} $w(0,y)=v(0,B)u(y) = u$, $\dom(v(x,B)B)\subset C^\infty_{loc} $. However, given $u\in C^\infty_0$, we do not know if $u\in \dom(v(x,B))$.
 \end{proof}
\section{Hypoellipticity and closed graph theorem}\label{sec:projections}

      % \hide{Ideally, we want to apply $v(x,B)$ to an arbitrary $u\in C^\infty_0$, but it may not be in the domain of $v(x,B)$.}
     For generic $u\in L^2$ or even $C^\infty_0$, we may not be able to guarantee the hypothesis of Lemma \ref{lem:tld-spectral}. To go around this concern, we first use the Closed Graph Theorem to control smoothness of solutions to $Lw=0$ using hypoellipticity. We start with a bit of notation.
     \begin{defn}\label{def-nesting}
       For $r>0$ define $I_r:=(x_0-r,x_0+r)$
     \end{defn}

     \begin{lem}\label{lem:X} Let $y_0\in R^n_y$ with $y_0\in \Omega \Subset \R^n_y$ and let $r>0$. For $I_r$ from Definition \ref{def-nesting}, define
     \begin{align}\label{hypo:space}
       X:=\{ w \in L^2({I_r}\times \Omega  )\mid L w = 0 \text{ in the weak sense on }I_{r}\times \Omega
        \}
     \end{align}
      Then $X$  is a closed subspace of $L^2_{loc}({I_r}\times \Omega  )$.
     \end{lem}
     \begin{proof}
       We invoke duality. Let $w_n \to w \in L^2_{loc}$ with $L w_n=0$ weakly. Then for $\phi \in C^\infty_0({I_r}\times \Omega  )$
  \begin{align*}
    (Lw,\phi) = (w, L^*\phi) = \lim_n (w_n, L^*\phi)=0
  \end{align*}
  Thus $Lw=0$ weakly and $w\in X$.
     \end{proof}

      \begin{lem}\label{closed-graph}[Closed Graph] Let $\Omega'\Subset\Omega \subset \R^n_y$ and $0<r'<r$. Let $s> 0$. Define $I_r$ and $I_{r'}$ by Definition \ref{def-nesting}.\\
  Suppose the operator $L$ from \eqref{tldL} is $H^s$ hypoelliptic at $(x_0,y_0)$ in the sense of Definition \ref{hypo-weak} and $X$ is as in \eqref{hypo:space}. Then the operator \begin{align*}
      \hspace{-20pt}  T:X \to H^s(I_{r'}\times \Omega') \text{ defined by } Tw = w
      \end{align*}
     has a closed graph. In particular, if $w \in X$, then
     \begin{align}
       \label{closed-graph:est} \norm{w}_{H^s(I_{r'}\times \Omega')}\le \tld C(\Omega,\Omega',r,r')\norm{w}_{L^2({I_{r}}\times \Omega)}
     \end{align}
  \end{lem}

  \begin{rem}\label{rem:closed-graph}
    As $r'\to r$ or $d(\partial\Omega',\partial\Omega)\to 0$ the constant in \eqref{closed-graph:est} may go to infinity.
  \end{rem}
  %Both the applications and the proof of the Lemma require $s>\frac{1}{2}$ to invoke Sobolev embedding for $x\in \R^1$.
  \begin{proof}
  First, we demonstrate $T$ is well defined. Let $\psi\in C^\infty_0(I_{r}\times \Omega)$ be a bump function $0\le \psi\le 1$ with $\psi(x,y)\equiv 1$ on $I_{r'}\times\Omega'$. Then, for all $u\in X$, $L(\psi u)(x,y)=0$ for $(x,y)\in I_{r'}\times \Omega'$. By Definition \ref{hypo-weak} $\psi u\in H^s$. Hence $u\in H^s(I_{r'}\times \Omega')$. \\

  \noindent To demonstrate that $T$ has a closed graph, suppose $w_n\in X$, $$w_n\to w\in X,  \text{ and }Tw_n\to \tld w \in H^s(I_{r'}\times \Omega'), $$
  Passing to subsequences if necessary, we may assume that $w_n\to \tld w$ and $w_n\to w$ for a.e. $(x,y)\in I_{r'}\times \Omega'$. Since $H^s(I_{r'}\times \Omega') \subset L^2(I_{r'}\times \Omega')$, it is then clear that $w=\tld w$ and the graph of $T$ is closed.
 By the closed graph theorem,   c.f. \cite{YosFunct} p.80, $T$ is continuous, and hence satisfies \eqref{closed-graph:est}.
  \end{proof}
     %As $H^s_{loc}$ is an F-space, rather than a Banach space,
We now need a version of Sobolev Embedding to apply Lemma \ref{closed-graph} in $H^s_y$ for $\frac{1}{2}<s\le\frac{n+1}{2}$, where $y\in \R^n$. If we use $C^\infty$ hypoellipticity or $H^s$ hypoellipticity for $s$ bigger than this range, the standard Sobolev Embedding $H^{\frac{n+1}{2}^+}(\R^{n+1})\subset C^0(\R^{n+1})$ suffices.% allows us to demonstrate sufficiency of \eqref{superlog} for our operator $L_2$.\\

      \begin{lem}\label{Sobolev-1D}[Sobolev Embedding]
        Let
        \begin{align}\label{s=s_1+}
          s=s_1+s_2\text{ for }s_1>\frac{1}{2}\text{ and }s_2\ge 0
        \end{align}
         Suppose $u(x,y)\in H^{s}_{x,y}(\R_x\times \R^n_y)$ . Then $x\to u(x,\cdot)$ is a continuous uniformly bounded function from $\R\to H^{s_2}_y$ with
      \begin{align}\label{Sobolev-estimate}
        \norm{u}_{L^\infty_xH^{s_2}_y}\le C_{s_1,s_2}\norm{u}_{H^s_{x,y}}
      \end{align}
        \end{lem}
        \begin{proof}
        Define
        \begin{align}\label{sobolev-u}
          u_s:=\jap{\dy}^{s_2} \jap{\dx}^{s_1}u,
        \end{align}
    where $\jap{\dx}^{s_1}$ is a pseudodifferential operator with symbol $\jap{\xi}^{s_1}$ and similarly for $\jap{\dy}^{s_2}$.
        We can re-write this formula as
          \begin{align*}
            \hat u_s(\xi,\eta)=\frac{(1+|\xi|^2)^{s_1}(1+|\eta|^2)^{s_2}}{(1+|\xi|^2+|\eta|^2))^{s}}(1+|\xi|^2+|\eta|^2))^{s}\hat u(\eta,\xi)
          \end{align*}
          By hypothesis, $u\in H^s$, and the term $\frac{(1+|\xi|^2)^{s_1}(1+|\eta|^2)^{s_2}}{(1+|\xi|^2+|\eta|^2))^{s}}$ is a bounded Fourier multiplier for $s_1$, $s_2\ge 0$. Thus
          \begin{align} \label{Sobolev-derivative}
            \norm{u_s}_{L^2_{x,y}}\le C_{s_1,s_2}\norm{u}_{H^s_{x,y}}
          \end{align}
          We return to the analysis of $u$ using \eqref{sobolev-u}. Write formally
        \begin{align}\label{Sobolev-u(x)}
           u(x,y)=\jap{\dx}^{-s_1}\jap{\dy}^{-s_2}u_s(x,y),
        \end{align}
          and define, $g(x)$ as an inverse Fourier transform of $\jap{\dx}^{-s_1}$, i.e.
          \begin{align*}
            \hat g(\xi)=(1+|\xi|^2)^{-s_1} \in L^2_\xi \text{ for }s_1>\frac{1}{2}.
          \end{align*}
          Note that from the definition $\norm{g}_{L^2_x}\le C_{s_1}$. With this notation we rewrite \eqref{Sobolev-u(x)} as
          \begin{align}\label{Sobolev-convolution}
           u(x,y) = g(x)\conv_x \jap{\dy}^{-s_2} u_s.
          \end{align}
          Applying Cauchy-Schwartz, \eqref{Sobolev-derivative} and the $L^2$ estimate for $g$ implies %
          \begin{align}\label{sobolev-cauchy-schwarz}
            \norm{u}_{L^\infty_xH^{s_2}_y} \le C_{s_1}\norm{\jap{\dy}^{-s_2} u_s}_{L^2_{x}H^{s_2}_y}=C_{s_1}\norm{u_s}_{L^2_{x,y}}\le C_{s_1,s_2}\norm{u}_{H^s_{x,y}},
          \end{align}
%           We conclude
which concludes \eqref{Sobolev-estimate}. Furthermore, to demonstrate continuity we use \eqref{Sobolev-convolution}
\begin{align*}%\label{Sobolev-difference}
\norm{u(x,\cdot)-u(x',\cdot)}_{H^{s_2}_y}\le \int \left| g(x-z) - g(x'-z)\right| \norm{\jap{\dy}^{-s_2} u_s(z,\cdot)}_{H^{s_2}_y} dz
\end{align*}
Then by Cauchy-Schwartz
\begin{align*}
  \norm{u(x,\cdot)-u(x',\cdot)}_{H^{s_2}_y} \le \left(\int|g(x-z) - g(x'-z)|^{2}dz\right)^{1/2}\norm{u_s}_{L^2_{z,y}}
\end{align*}
As $x'\to x$, the first factor on the right converges to zero by properties of the translation operator in $L^2$.\qedhere
%\begin{align*}
%  G(x,x',z)= \left| g(x-z) - g(x'-z)\right| \norm{\jap{\dy}^{-s_2} u_s(z,\cdot)}_{H^{s_2}_y}
%\end{align*}
%Note that $G\in L^1_z$ by an argument identical to \eqref{sobolev-cauchy-schwarz} and $\lim_{x'\to x}G(x',x,z)=0$ pointwise in $z$. By the Lebesgue Dominated convergence, $u(x,\cdot)-u(x',\cdot)\to 0$ in $H^{s_2}_y$.
%As
%\begin{align*}
%  \int_{|z|>R} G(x,x',z) dz\to 0 \
%\end{align*}

%\[
%\norm{u(x,\cdot)-u(x',\cdot)}_{H^{s_2}_y}\le 2C_{s_1,s_2}\norm{u}_{H^s_{x,y}}\  \  \forall x,x'\in \mathbb{R},
%\]
%{\LARGE but we want to show this quantity goes to zero as $x\to x'$. I don't see how to use the Lebesgue Dominated convergence here}
%%and since $u\in H^{s}_{x,y}(\R_x\times \R^n_y)$,
%%\[
%%
%%\]
%
%the Lebesgue Dominated Convergence Theorem for $\norm{u(x,\cdot)-u(x',\cdot)}_{H^{s_2}_y}$, the map $x\to u(x,\cdot)$ is uniformly continuous from $\R_x \to H^{s_2}_y$.
        \end{proof}
        We conclude the section by combining Lemmas \ref{closed-graph} and \ref{Sobolev-1D}.
        \begin{prop}\label{prop:hypo-closed}
          Let $s=s_1+s_2$ with $s_1>\frac{1}{2}$ and $s_2> 0$. Suppose the operator $L$ from \eqref{tldL} is $H^s$ hypoelliptic at $(x_0,y_0)$. Then for some $\Omega'\Subset\Omega\Subset \R^n_y$ and $r>0$ small enough,
             \begin{align}
       \label{hypo:est} \norm{w(x_0,\cdot)}_{H^{s_2}(\Omega')}\le \tld C(s_1,s_2,\Omega,\Omega',r)\norm{w}_{L^2({I_{r}}\times \Omega)}
     \end{align}
        \end{prop}
        \section{Hypoellipticity with spectral projections}\label{projection:gain}
     % {\bf I wrote a sketch of the proof, but need a bit more time to tidy it up. I certainly want to write it As we discussed today, by raising $s$ from $\frac{1}{2}^+$ to $\frac{n+1}{2}^+$ we can avoid anisotropic Sobolev embedding. I will clean up the Closed Graph statement as well}
      As mentioned above, for generic $u\in L^2_y$ or even $C^\infty_0$, we may not be able to guarantee the hypothesis of Lemma \ref{lem:tld-spectral}. We show below that for hypoelliptic operators, Lemma \ref{lem:tld-spectral} applies for any $P_ju$. In particular, this will imply that for hypoelliptic operator $L$ the spectral projection $P_j$ from \eqref{P-j-def} is smoothing. We state this in a form that is crucial for the proof of Theorem \ref{thm:main}.
\begin{prop}\label{lem:necess:Hs}
Let $s=s_1+s_2$ with $s_1>\frac{1}{2}$ and $s_2> 0$. Suppose $L$ from  \eqref{tldL} is $H^s$ hypoelliptic near $(x_0,y_0)$ and let the operator $P_j$ be defined by \eqref{P-j-def}. Let $\eps>0$ and $y_0\in \Omega'\Subset\Omega\Subset \R^n_y$, then
       \begin{align}\label{Hs-tldL}
       \norm{P_j u}_{H^{s_2}(\Omega')}  \le \tld C(\eps,x_0,s_1,s_2,\Omega,\Omega')  \exp(\eps \sqrt{e^j}) \norm{P_j u}_{L^2(\Omega)}
     \end{align}
     for all $u(y)\in L^2_y(\Omega)$ and all $j\ge 0$. The constant $\tld C\to \infty$ as $\eps\to 0$, and similarly as $d(\partial\Omega',\partial\Omega)\to 0$, see Remark \ref{rem:closed-graph}.
\end{prop}
%\hide{Essentially the log estimate is the interpolation between \eqref{necess:Hs} for high (Fourier) frequencies and the operator $A$ for small frequencies. Being able to split frequencies is where Littlewood-Paley works well}
% We now want to obtain \eqref{superlog} for $(L_2 u ,u) $ from $v(x,B)u$.\\
 \begin{proof}
Let $u\in L^2$ and define $u_j=P_ju$ and $$w_j = v(x,B)u_j= \psi_j(B) v(x,B) u.$$ Note by \eqref{lambda-ode} $w_j(x_0)=u_j$. By Lemma \ref{lem:LP-operator}
 \begin{align}\label{Pj-exp}
   \norm{ e^{\eps \sqrt{B}} u_j(y) }_{L^2_y} \le \sum_{|j'-j|\le 1} e^{\eps \sqrt{e^{j'}}} \norm{P_{j'} P_j u}<\infty
%   \le e^{\eps \sqrt{e^{j+1}}}\norm{P_j u}
 \end{align}
Therefore, Lemma \ref{lem:tld-spectral} applies to $u_j$ and gives $$L w_j(x,y)=0 \text{ for } |x-x_0|\le r_0.$$ Thus hypoellipticity of $L$ implies that $w_j \in C^\infty_{loc}$. Note, in particular, $P_ju$ is locally smooth. It remains to establish \eqref{Hs-tldL} - a quantitative control of this smoothness.\\

We now invoke Proposition \ref{prop:hypo-closed}. By \eqref{hypo:est}
 \begin{align}\label{Hs:to:L2}
   \norm{w_j(x_0,\cdot)}_{H^{s_2}(\Omega')} \le \tld C \norm{w_j}_{L^2(I_r\times \Omega)}
 \end{align}
where the constant $\tld C$ is independent of $w_j$. In particular, as $r\to 0$, $\tld C$ may go to infinity by Remark \ref{rem:closed-graph}.\\

Recall that $w_j(x_0,y)=P_ju$. Therefore we can restate \eqref{Hs:to:L2} as
\begin{align}
  \label{closed:Hs}
   \norm{P_ju}_{H^{s_2}(\Omega')} \le \tld C \norm{w_j}_{L^2(I_r\times \Omega)}
\end{align}
 %for all $w(x,y)\in \Schw'$.
Proposition \ref{spectral} and the support of $\psi_j$ imply
\begin{align}\label{Hs-quant}
  \norm{w_j(x,y)}_{L^2(I_r\times \Omega)}^2 \le\int_{x_0-r}^{x_0+r}\int_{e^{j-1}}^{e^{j+1}} |v(x,\lambda)|^2 \psi_j^2(\lambda) d\norm{E_\lambda u}_{L^2_y}^2dx%\les_{x_0}% \exp(C_{x_0}|x-x_0| \sqrt{e^j})\norm{P_ju}^2
\end{align}
From Lemma \ref{lem:ODE} we get
\begin{align*}
  |v(x,\lambda)| \le e^{C_{x_0}\jap{\lambda}^{\frac{1}{2}}|x-x_0|}
\end{align*}
Note, that $C_{x_0}$, unlike $\tld C$ is independent of $\eps$. Hence for $\lambda \in [{e^{j-1}},{e^{j+1}}]$ and $x\in [{x_0-r},{x_0+r}]$:
\begin{align*}
  & |v(x,\lambda)|^2\le e^{2 C_{x_0}e^{\frac{j+1}{2}}r}=e^{ C_{x_0}e^{\frac{j}{2}}r}, \text{ where } C_{x_0} \text{ is replaced with a larger constant}\\
  & \text{that is still independent of $j$ and $\eps$ in the last equality.}
\end{align*}
%We can relabel the constant $C_{x_0}$ to turn the estimate into $|v(x,\lambda)|^2\le e^{ C_{x_0}e^{\frac{j}{2}}r}$.
 For $\eps>0$ in the hypothesis of our Proposition choose
\begin{align*}
  r = C_{x_0}^{-1}\eps.
\end{align*}
Then
\begin{align*}
  |v(x,\lambda)|^2 \le  \exp(\eps \sqrt{e^j}).
\end{align*}
We now return to \eqref{Hs-quant} with this estimate of $v$ on our projected band:
\begin{align*}
 & \norm{w_j(x,y)}_{L^2(I_r\times \Omega)}^2 \le \int_{x_0-r}^{x_0+r} \exp(\eps \sqrt{e^j})\int_{e^{j-1}}^{e^{j+1}} \psi_j^2(\lambda) d\norm{E_\lambda u}_{L^2_y}^2dx
\end{align*}
Since $ \norm{P_ju}_{L^2(\Omega)}^2= \int_{e^{j-1}}^{e^{j+1}} \psi_j^2(\lambda) d\norm{E_\lambda u}_{L^2(\Omega)}^2$, we rewrite this as
\begin{align*}
  \norm{w_j}_{L^2(I_r\times \Omega)}^2 \le  \eps \exp(\eps \sqrt{e^j}) \norm{P_ju}^2_{L^2(\Omega)}
\end{align*}
                                        % \begin{align*}%\label{w_j}
%   \norm{w_j(x,y)}_{L^2_{x,y,loc}} \le C(x_0)\exp(\eps \sqrt{e^j})\norm{P_ju}
%                                                       \end{align*}
This estimate together with \eqref{closed:Hs} gives
\begin{align*}
  \norm{P_ju}_{H^{s_2}(\Omega')}^2 \le \tld C  \eps \exp(\eps \sqrt{e^j}) \norm{P_ju}^2_{L^2(\Omega)}
\end{align*}
Note, that this estimate does not imply that $P_ju\to 0$, as $\eps\to 0$, since $\tld C$ depends on $\eps$ itself. Absorbing, $\eps$ into a new constant $\tld C$ concludes the proof.
  \end{proof}
  \section{Interpolation and proof Theorem \ref{thm:main}}\label{sec:interpolation}
\begin{lem}  \label{lem:interp}
    Let $\eps>0$, $s_2>0$. For each $\xi\in \R^n$ let $R=R(\xi)$ be defined by
    \begin{align}\label{R-xi}
      e^R = \left(\frac{s_2\log\jxi}{2\eps}\right)^2.
    \end{align}
   Consider a sequence of functions $\alpha_k(y) \in H^{s_2}_y$ for $k=0, 1, 2 \ldots$ Then
    \begin{align}
      \label{low-j}& \norm{\sum_{k\le R(\xi)} \log\jxi \widehat{\alpha_k}(\xi)}_{L^2_\xi}^2 \le C(\eps,s_2)\sum_{k=0}^\infty e^{-2\eps \sqrt{e^k}}\norm{\alpha_k}_{H^{s_2}_y}^2\\
      \label{high-j}& \norm{\sum_{k\ge R} \log\jxi \widehat{\alpha_k}(\xi)}_{L^2_\xi}^2 \le C\left(\frac{\eps}{s_2}\right)^2\sum_{k=0}^\infty e^k\norm{\alpha_k}_{L^2_y}^2
    \end{align}
  \end{lem}
  \begin{proof}
  Both of the inequalities are proved by a clever use of the Cauchy-Schwartz inequality.\\

   By the Cauchy-Schwartz inequality
    \begin{align*}
      \sum_{k\le R} \log\jxi| \widehat{\alpha_k}(\xi)| \le \left(\sum_{k\le R} \frac{\log^2\jxi}{\jxi^{2s_2}} e^{2\eps \sqrt{e^k}}\right)^{\frac{1}{2}}\cdot \left( \sum_{k\le R}\frac{\jxi^{2s_2}}{e^{2\eps \sqrt{e^k}}} |\widehat{\alpha_k}|^2\right)^{\frac{1}{2}}
    \end{align*}
    The terms of the first sum are increasing in $k$, so the sum can be estimated
    \begin{align*}
      \sum_{k\le R} \frac{\log^2\jxi}{\jxi^{2s_2}} e^{2\eps \sqrt{e^k}} \le \frac{\log^2\jxi}{\jxi^{2s_2}} R \cdot e^{2\eps \sqrt{e^R}}
    \end{align*}
     Substitution of the identity \eqref{R-xi} (which is equivalent to $\jxi^{s_2} = e^{2\eps \sqrt{e^R}}$) implies
     \begin{align*}
       \sum_{k\le R} \frac{\log^2\jxi}{\jxi^{2s_2}} e^{2\eps \sqrt{e^k}} \le \frac{\log^2\jxi (2\log s_2+2\log \log \jxi-2\log2\eps)}{\jxi^{s_2}} \le C(\eps,s_2)
     \end{align*}
 Combining the last three estimates we obtain
    \begin{align*}%\label{low-j-split}
      \norm{\sum_{k\le R} \log\jxi \widehat{\alpha_k}(\xi)}_{L^2_\xi}^2 \le  C(\eps,s_2) \int \sum_{k\le R}\frac{\jxi^{2s_2}}{e^{2\eps \sqrt{e^k}}} |\widehat{\alpha_k}|^2 d\xi
      % \le \norm{\sum_{k\le R} \frac{\log^2\jxi}{\jxi^{2s_2}} e^{2\eps \sqrt{e^k}}}_{L^\infty_\xi}\cdot \sum_{k\le R}\int\frac{\jxi^{2s_2}}{e^{2\eps \sqrt{e^k}}} |\widehat{\alpha_k}|^2d\xi
    \end{align*}
  Observe
     \begin{align*}
      \sum_{k\le R}\frac{\jxi^{2s_2}}{e^{2\eps \sqrt{e^k}}}|\widehat{\alpha_k}|^2 \le  \sum_{k=0}^\infty\frac{\jxi^{2s_2}}{e^{2\eps \sqrt{e^k}}}|\widehat{\alpha_k}|^2
     \end{align*}
Interchanging the sum and integration using Fubini concludes the proof of \eqref{low-j}.

  Similarly for \eqref{high-j}, by Cauchy-Schwartz
  \begin{align*}
    \left(\sum_{k\ge R} \log\jxi| \widehat{\alpha_k}(\xi)| \right)^2\le \left( \sum_{k\ge R} \frac{\log^2\jxi}{e^k}\right) \cdot \left( \sum_{k\ge R} e^k |\widehat{\alpha_k}|^2\right)
  \end{align*}
  The first sum is a geometric series in $k$, so
  \begin{align*}
    \sum_{k\ge R} \frac{\log^2\jxi}{e^k} = \frac{\log^2\jxi}{e^R}\frac{e}{e-1} \le 3\left(4\frac{\eps}{s_2}\right)^2
  \end{align*}
Integrating in $\xi$ gives
   \begin{align*}
    \norm{\sum_{k\ge R} \log\jxi \widehat{\alpha_k}}_{L^2}^2 \le C\left(\frac{\eps}{s_2}\right)^2\int \sum_k e^k|\hat \alpha_k|^2 d\xi
  \end{align*}
 Applying Fubini concludes \eqref{high-j}
  \end{proof}
  %We are now ready to conclude the Theorem \ref{thm:main}.
 \begin{proof}[Proof of Theorem \ref{thm:main}]
 Suppose $L$ from \eqref{tldL} is $H^s$ hypoelliptic for $s>\frac{1}{2}$, then Proposition \ref{lem:necess:Hs} applies for $s_2=s-s_1>0$ with some $s_1>\frac{1}{2}$.\\

   Let $u(y)\in C^\infty_0(\Omega')$, and define $\alpha_j =  P_ju$. Since $P_j$'s are a resolution of identity by Lemma \ref{lem:LP-operator}, we can decompose
   \begin{align*}
   u=  \sum_{j\ge 0} P_j u= \sum_{j\ge 0} \alpha_j
   \end{align*}
  By the triangle inequality $$ \norm{\log\jxi\hat u}_{L^2} \le  \norm{\log\jxi \sum_{j\le R}\hat \alpha_j}_{L^2}+  \norm{\log\jxi \sum_{j\ge R} \hat \alpha_j}_{L^2}$$
  By Cauchy-Schwartz
  \begin{align*}
    \norm{\log\jxi \hat u}_{L^2}^2\le 2\left(\norm{\log\jxi \sum_{j\le R} \hat\alpha_j }_{L^2}^2 + \norm{\log\jxi \sum_{j\ge R}\hat \alpha_j }_{L^2}^2\right)
  \end{align*}
  We are now set to use Lemma \ref{lem:interp} for our choice of $\alpha_j$ (replacing $k$ with $j$):
  \begin{align}\label{log-necess-split}
     \norm{\log\jxi u}_{L^2}^2  \le C\frac{\eps}{s_2}^2\sum_{j=0}^\infty e^j\norm{ P_ju }_{L_2}^2 + C_\eps\sum_{j=0}^\infty e^{-2\eps \sqrt{e^j}}\norm{ P_j u}_{H^{s_2}}^2 := I+ II
  \end{align}
  We estimate $I$ and $II$ separately. For $I$ we estimate the first sum via \eqref{LP-oper} for $f(\lambda)=\sqrt{\lambda}$. As $f(\lambda)>0$ and increasing for $\lambda\ge 1$
\begin{align*}
 I= C\eps^2\sum_{j\ge R} e\cdot |\sqrt{e^{j-1}}|^2\norm{P_ju }_{L_2}^2 \le Ce\eps^2 \norm{\sqrt{B} u}^2
\end{align*}
Choose $\eps>0$ small enough, so that $Ce\eps\le 1$. Then
\begin{align*}
  I\le \eps \norm{\sqrt{B} u}^2
\end{align*}
Using $(Bu,u) = (L_2 u, u)$ and the fact that $B$ is an extension of the operator $L_2$ and $u\in C^\infty_0$ we obtain
\begin{align*}
  I\le \eps (L_2u,u)
\end{align*}

For $II$ we %first apply \eqref{low-j} from Lemma \ref{lem:interp}:
%\begin{align*}
%  II\le   C_\eps\sum_{j=0}^\infty e^{-2\eps \sqrt{e^j}}\norm{\alpha_j}_{H^{s_2}_y}^2
%\end{align*}
 observe from Proposition \ref{lem:necess:Hs}
\begin{align*}
  \norm{P_j u}_{H^{s_2}}^2 \le C(\eps,s) e^{2\eps \sqrt{e^j}}\norm{P_j u}_{L^2}^2
\end{align*}
Hence
\begin{align*}
  II \le C(\eps,s) \sum_j \norm{ P_j u}_{L^2}^2 \le C_\eps \norm{u}^2
\end{align*}
Combining $I$ and $II$ into \eqref{log-necess-split} we obtain
\begin{align*}
  \norm{\log\jxi u}_{L^2(\Omega')}^2 \le \eps (L_2 u, u) + C(\eps,s) \norm{u}^2_{L^2(\Omega)}\hspace{20pt}\qedhere
\end{align*}
 \end{proof}
\begin{appendices}

\section{Spectral projections for the operator}\label{sec:spectral}
We start with recalling the relevant facts about the spectral projection and constructing Littlewood-Paley projections adapted to the operator.\\

%\begin{enumerate}
Note that by \eqref{L2}, $L_2$ is an operator bounded below. That is, given $\Omega \Subset \R^n$ there is a constant $c_\Omega\in \R$ so that
 \begin{align}\label{global-lower-bound}
  (L_2 u, u) \ge c_{\Omega}\norm{u}^2 \text{ for all } u\in C_0^\infty(\Omega)
\end{align}

  Fix a neighborhood of $(x_0,y_0)$ with compact closure. Since Theorem \ref{thm:main} is a local result, we may change the operator $L$ from \eqref{tldL} away from this neighborhood to assume \eqref{global-lower-bound} holds with a constant $c=c_{\Omega_{(x_0,y_0)}}$.
\begin{rem}\label{rem:lambda-one}
  We may assume without loss of generality that \eqref{global-lower-bound} holds for $c=1$ and hence $L_2$ is a positive operator:
  \begin{align}\label{L2:positive}
 (L_2u,u) \ge \norm{u}^2
\end{align}
\end{rem}
\begin{proof}
  If $c<1$ in \eqref{global-lower-bound}, we replace $L_2$ with $L_2:=L_2-c+1$ for $c$ from \eqref{global-lower-bound} and $a_0(x)$ with $a_0(x)+g(x)[c-1]$. Note that with this change the operator $L$ in \eqref{tldL} remains unchanged. Such a replacement makes $L_2$ a positive operator, which we would still call $L_2$.
\end{proof}

%This has no effect of $L$ if we then replace $a_0(x)$ with $a_0(x)-c$.

 We now extend the domain of $L_2$ from $C^\infty_0(\Omega)$ to a unique maximal domain inside $L^2(\Omega)$ to make it self adjoint. More precisely, let $B$ be the Friedreich extension of $L_2$ (c.f. Theorem 2 on p. 317 of \cite{YosFunct}). We then have that $B$ satisfies
  \begin{align}\label{positivity-operator}
   (Bu,u) \ge \norm{u}^2,
 \end{align}
since $B u = L_2 u$ for $u\in C^\infty_0$. Moreover, as $B$ is a self-adjoint operator, the spectral theorem/calculus applies to it. % In particular, $B$ has a functional calculus defined.
  We summarize the relevant parts from \cite{YosFunct}  ch XI sec. 5 and 10 (p.300-343)]\

 % I.e. there exists a Banach space  $Bu = L_2 u$ for all $u\in C^\infty_0$ and $B$ is self-adjoint on $dom B$. \hide{basically enlarge domain from $C_0^\infty$ to a maximal subset of $L_2$ on which it is self-adjoint)}. %Then take a  square root of that operator. \cite{YosFunct} Ch. XI.
    %\hide{Without completing the domain of the operator, we cannot do spectral projection}% arguments, which are often relevant}
  % Since $B$ is self-adjoint, it has a spectral resolution with the following properties:

   \begin{prop}\label{spectral}
   There exists a spectral projection $E_B(\lambda):L^2(\Omega)\to L^2(\Omega)$, that is a resolution of identity and allows to form functions of the operator $B$ as follows:
   \begin{enumerate}
  \item \label{proj} $E_B$ is a projection, i.e. \begin{align}\label{spectral-resolution}
    E_B(\lambda) E_B(\lambda') = E_B\left[\min(\lambda,\lambda')\right]
  \end{align}
   And a resolution of identity: $E_B(\infty) = I$, $E_B(-\infty)=0$; Moreover $E$ is right-continuous.
   \item \label{spectral-norm} For a measurable function $f(x)$ we can define the function $f(B)$ of the operator $B$ as follows:
\begin{align}\label{spectral-function}
f(B)u = \int_{-\infty}^\infty f(\lambda) dE_B(\lambda) u
\end{align}
Where $\dom f(B) = \{ u\in L^2(\Omega): \norm{f(B)u}^2:= \int_{-\infty}^\infty |f(\lambda)|^2 d ||E_B(\lambda) u||^2 <\infty \}.$
   \item In particular, for $f(\lambda)=1$ and $\lambda$ gives the following identities
    \begin{align}\label{spectral-identity}
     u = \int_{-\infty}^\infty d E_B(\lambda) u & \hspace{10pt} B u=\int_{-\infty}^\infty \lambda d E_B(\lambda) u
   \end{align}
   are defined on $L^2(\Omega)$ and $\dom(B)$ respectively.
\item \label{spectral:commutes} For any measurable functions $f$ and $g:\R\to \R$, $f(B)$ and $g(B)$ commute on the $dom \left(f(B) g(B)\right)$ defined via \eqref{spectral-function}. Moreover,
\begin{align} \label{spectral:commutes:eq}
  (f(B)g(B)u,v)=(f(B)u,g(B)v) = \int f(\lambda)\overline{g(\lambda)} d(E_B(\lambda) u,v)
\end{align}
  \end{enumerate}

\end{prop}
\begin{proof}
  Existence of spectral projection is given Theorem XI.6.1 on p. 313 of \cite{YosFunct} via unitary operators. Part 2 and Part 3 are proved in Theorem XI.5.2. on p.311 using Part 1 and orthogonality. Part 4 is proved in Corollary XI.5.2 and for measurable functions in Theorem XI.12.3 on p. 343 of \cite{YosFunct}.
\end{proof}
We follow the notation of \cite{YosFunct} and denote
\begin{align}\label{proj-band}
  E(\lambda,\mu):=E(\lambda)-E(\mu)\text{ for }\lambda \ge \mu
\end{align}

One of the important consequences of part 4 of Proposition \ref{spectral} is the following equality $$\norm{\sqrt{B}u}^2 = (Bu,u).$$
\begin{lem}\label{lem:positive}
  Suppose \eqref{positivity-operator} holds for $B$. Then the projection $E_B$ vanishes below $\lambda=1$, i.e. $E_B(1^-)=0$,
\end{lem}
\begin{proof}
It suffices to show that $u:=E(\mu)v=0$ for all $\mu <1$ and all $v\in L^2(\Omega)$. To achieve that we substitute $u$ into \eqref{positivity-operator} and \eqref{spectral-identity} to obtain
$$
0\le \int_{-\infty}^\infty (\lambda-1) d (E_B(\lambda) u, u)
$$
We first consider the large frequencies and decompose the integral into a Riemann sum as follows
\begin{align*}
  & \int^{\infty}_\mu (\lambda-1) d(E(\lambda) u,  u)\\
  & = \sup_{\mu=\lambda_0<\lambda_1<\ldots<\lambda_N=M}\sum_{j=0}^N (\lambda_j-1) (E(\lambda_j,\lambda_{j-1}) E(\mu)v, E(\mu)v),
\end{align*}
Here the supremum is taken over all partitions $\mu=\lambda_0<\lambda_1<\ldots<\lambda_N=M$ of all intervals $[\mu,M]\subset [\mu,\infty)$. By \eqref{proj-band} and \eqref{spectral-resolution} all terms in the sum above vanish. We conclude that
\begin{align*}
   & \int^{\infty}_\mu (\lambda-1) d(E(\lambda) u,  u) = 0
\end{align*}

Hence the integral over the low frequencies must be non-negative. We proceed in the low frequency case as before
\begin{align*}
\int_{-\infty}^\mu (\lambda-1) d(E(\lambda) u,  u) = \sup_{-M=\lambda_0<\lambda_1<\ldots<\lambda_N=\mu}\sum_{j=1}^N (\lambda_j-1) (E(\lambda_j,\lambda_{j-1}) E(\mu)v, E(\mu)v)
\end{align*}
with the supremum taken over all partitions $-M=\lambda_0<\lambda_1<\ldots<\lambda_N=\mu$ of all intervals $[-M,\mu]\subset (-\infty,\mu]$. Observe that \eqref{spectral-resolution} implies orthogonality of projections localized at different frequencies, i.e. $E(\lambda_j,\lambda_{j-1})E(\lambda_k,\lambda_{k-1})u=0$ for different $j$ and $k$. Moreover, using (\ref{spectral:commutes:eq}) with $f=\chi_{(0,\lambda_j)}-\chi_{(0,\lambda_{j-1})}$ and $g=\chi_{(0,\lambda_k)}-\chi_{(0,\lambda_{k-1})}$ we have
\begin{align*}
\int_{-\infty}^\mu (\lambda-1) d(E(\lambda) u,  u)= \sup_{-M=\lambda_0<\lambda_1<\ldots<\lambda_N=\mu}\sum_{j=0}^N (\lambda_j-1) \norm{E(\lambda_j,\lambda_{j-1}) v}^2\\
 \le (\mu-1)\sup_{-M=\lambda_0<\lambda_1<\ldots<\lambda_N=\mu}\sum_{j=0}^N \norm{E(\lambda_j,\lambda_{j-1}) v}^2
 %\text{ which leads to a contradiction for }\mu<1.
\end{align*}
Using orthogonality of distinct frequency interval again we obtain
\begin{align*}
   \int_{-\infty}^\mu (\lambda-1) d(E(\lambda) u,  u) \le (\mu-1)\norm{E(\mu) v}^2
\end{align*}
Combining all the estimates established in this proof we obtain
\begin{align*}
  0\le \int_{-\infty}^\infty (\lambda-1) d (E_B(\lambda) u, u) \le (\mu-1)\norm{E(\mu) v}^2\le 0,\quad\text{since}\  \  \mu<1.
\end{align*}
We must conclude that $E(\mu) v=0$.
\end{proof}

     % This operation is well defined for every self-adjoint operator. Note that $E_B(\lambda)$ commutes with $B$ for every $\lambda$.
In light of Lemma \ref{lem:positive}, we only need to consider $\lambda\ge 1$.
For the analysis below it is convenient to localize solutions to specific frequencies. Littlewood-Paley decomposition is ideal for that. Let $\phi$ be a fixed cut-off function satisfying
\begin{align*}
  \phi(\lambda)\in C^\infty_0([-e,e]), \,\, \phi(\lambda)\equiv 1 \text{ for }|\lambda|\le 1 \text{ and } 0\le \phi \le 1.
\end{align*}
Moreover, we assume that $\phi$ is even and nonincreasing in $[0,e]$. Now, for any integer $j$ let
 \begin{align*}
   \phi_j(\lambda)\equiv \phi(\lambda\cdot e^{-j}).
 \end{align*}
Then $\phi_j(\lambda)=1$ for $|\lambda|\le e^j$ and $\phi_j(\lambda)=0$ for $|\lambda|\geq e^{j+1}$. We further define cut-offs localized to $|\lambda|\approx e^j$ as follows:
      \begin{align}\label{psi-j}
      \psi_j(\lambda) = \phi_j(\lambda)-\phi_{j-1}(\lambda) \text{ for } j\ge 0,
      \end{align}
 so that $\supp{\psi_j(|\lambda|)}=[e^{j-1},e^{j+1}]$ and $0\le \psi_j\le 1$.\\
         We summarize properties of such cutoffs in this range of $\lambda\ge 1$ as follows.
  \begin{lem}\label{Littlewood}
  Let $\psi_j$ be as defined in \eqref{psi-j} and consider $\lambda\ge 1$. We then have
  \begin{align}\label{LP-split}
  & \supp \psi_j(|\lambda|) = [e^{j-1},e^{j+1}] \text{ for } j\ge 0 ;   \hspace{10pt} \sum_{j=0}^\infty \psi_j(\lambda) = 1 \\\label{LP-ortho}
   & \psi_j(\lambda)\cdot \psi_k(\lambda) >0 \text{ if and only if } |j-k|\le 1; \hspace{10pt}   \frac{1}{2}\le \sum_{j=0}^\infty\psi_j^2 \le  1\
   \end{align}
   Furthermore, for any positive increasing function $f:[1,\infty)\to \R^{+}$
   \begin{align}
    \label{LP-norm}&
   \sum_{j=0}^\infty f(e^{j-1}) \psi_j^2(\lambda) \le f(\lambda)  \le 2\sum_{j=0}^\infty f(e^{j+1}) \psi_j^2(\lambda)
  % \\
  %&  \supp \psi_0(|\lambda|)=[0,e]\text{, and }\,\,  f(1)\psi_0^2(\lambda) \le  f(\lambda) \psi_0^2(\lambda) \le f(e) \psi_0^2(\lambda) \label{LP-0}
  \end{align}

   \end{lem}
\begin{proof}
From the definition of $\psi_j$ and supports of $\phi_j(\lambda)$, we immediately get the support property \eqref{LP-split} from  \eqref{psi-j}. %Moreover, for $\lambda\le e^{j-2}$, both terms in \eqref{psi-j} are equal to 1, so $\psi_j(\lambda)=0$.
%
% A computation shows for $j\ge 1$
% \begin{align}\bigotimes
%    \psi_j(e^{j-1})  = \phi(1/e)- \phi(1) = 0\\
%    \psi_j(e^j) = \phi(1) -\phi(e) = 1\\
%   \psi_j(e^{j+1})= \phi(e) - \phi(e^2) =0
%   \end{align}
%   and same for $\lambda$ in between.
A telescoping sum argument establishes the series identity in \eqref{LP-split}. Indeed
   \begin{align}
   \sum_{j=0}^n \psi_j (\lambda) = \phi_n(\lambda)- \phi_{-1}(\lambda)
 \end{align}
Since $\phi_n(\lambda)=1$ for $|\lambda|\le e^{n}$ and $\phi_{-1}(\lambda)= 0$ for $\lambda\ge 1$ , \eqref{LP-split} follows. %The rest is immediate from the overlap of supports and monotonicity of $f$.\\

 For \eqref{LP-ortho}, the orthonormal property is immediate from the fact that $\psi_j(\lambda)$ is nonzero only when $|\lambda|\in [e^{j-1},e^{j+1}]$. From the size of $\psi_j$, $\sum_{j=0}^\infty \psi_j^2 \le \sum_{j=0}^\infty \psi_j=1$. By (\ref{LP-ortho}) and the Cauchy-Schwartz inequality,
 $$
1=(\sum_{j=0}^\infty \psi_j(\lambda) )^2\le \sum_{j=0}^\infty \sum_{k=j-1}^{j+1} \psi_j\psi_k\le \sum_{j=0}^\infty( \psi_j^2 + \frac{1}{2}\psi_{j-1}^2+\frac{1}{2}\psi_{j+1}^2)
 $$
 which concludes \eqref{LP-ortho}.\\

 On the support of each $\psi_j$ for $j\ge 1$,  monotonicity of $f$ implies $$f(e^{j-1})\le f(\lambda)\le f(e^{j+1}).$$ Multiplying by $\psi_j(\lambda)^2$ makes the inequality valid for all $\lambda$:
$$ f(e^{j-1}) \psi_j^2(\lambda) \le  f(\lambda) \psi_j^2 \le f(e^{j+1}) \psi_j^2(\lambda).$$
 Summing up gives
 \begin{align*}
\sum_{j=0}^\infty f(e^{j-1}) \psi_j^2(\lambda) \le  \sum_{j=0}^\infty f(\lambda) \psi_j^2 \le \sum_{j=0}^\infty f(e^{j+1}) \psi_j^2(\lambda)
 \end{align*}
Finally, from \eqref{LP-ortho}
\begin{align*}
  \sum_{j=0}^\infty f(e^{j-1}) \psi_j^2(\lambda)  \le f(\lambda) \le 2 \sum_{j=0}^\infty f(e^{j+1}) \psi_j^2(\lambda)
\end{align*}

%
%\eqref{LP-0} is established similarly, with the hypothesis of $\lambda\ge 1$.
% The estimate $0<f(0)\approx f(1)<\infty$ concludes the $j=0$ case.
\end{proof}
With $\psi_j$ as above, we define a frequency localization adapted to the operator via Proposition \ref{spectral} \begin{align}\label{P-j-def}
        P_j u := \psi_j(B)u= \int \psi_j(\lambda) dE_B(\lambda)u %= \int_{e^{j-1}}^{e^{j+1}} \psi_j(\lambda) dE_B(\lambda)u
      \end{align} for $j\in \N$. Combining Proposition \ref{spectral}, Lemma \ref{lem:positive} and Lemma \ref{Littlewood} we get
   \begin{lem}  \label{lem:LP-operator}  Let $P_j$ be defined by \eqref{P-j-def}. Then for any $u\in L^2(\Omega)$ there holds
      \begin{align*}
        & u=\sum_{j=0}^\infty P_j u \hspace{15 pt} \norm{u}^2\approx \sum_{j=0}^\infty \norm{P_ju}^2;\\
        & \int P_j u P_{j'} v dx = 0 \text{ for } |j-j'|>1 \text{ and } v\in L^2
      \end{align*}
      %   $$ \hspace{15pt} $$% and this gives us an analogue of Littlewood-Paley decomposition adapted to our operator. Note, that for $j\ge 1$, integrating by parts
   Finally, for $f:[1,\infty)\to \mathbb{R}_+$ nondecreasing
\begin{align}\label{LP-oper}
  \sum_{j=0}^\infty |f\left(e^{j- 1}\right)|^2 \norm{P_j u}^2 \le \norm{f(B) u}^2 \le 2 \sum_{j=0}^\infty |f\left(e^{j+ 1}\right)|^2 \norm{P_j u}^2
\end{align}
%and the norm is finite if and only if the right side of the equation is finite.
      \end{lem}
      \begin{proof}

       By Proposition \ref{spectral} Parts 2 and 4, Lemma \ref{lem:positive} and \eqref{LP-split}
       \begin{align*}
         u = \int_{-\infty}^\infty dE(\lambda) u =\int_{1}^\infty \sum_{j=0}^\infty \psi_j(\lambda) dE(\lambda) u = \sum_{j=0}^\infty P_j u
                \end{align*}
                Similarly, using \eqref{LP-ortho}
       \begin{align*}
         \norm{u}^2 = \int_{-\infty}^\infty (dE(\lambda) u,u) \approx \int_{1}^\infty \sum_{j=0}^\infty  \psi_j^2(\lambda)(dE(\lambda) u,u) =  \sum_{j=0}^\infty \norm{P_ju}^2
       \end{align*}
       Orthogonality part of \eqref{LP-ortho} for $|j-j'|>1$ and \eqref{spectral:commutes:eq} establish
       $$\int P_j u P_{j'} v dx = 0$$
       Finally, from \eqref{LP-norm} and Fubini we conclude
       \begin{align*}
          \norm{f({B}) u}^2 = \int_{1}^\infty |f(\lambda)|^2 d (E_\lambda u,u) \le 2\sum_{j=0}^\infty \int_{1}^\infty f(e^{j+ 1})^2 \psi_j^2(\lambda)d (E_\lambda u,u)
       \end{align*}
       The lower bound is established similarly, completing the proof.
       % Then using \eqref{LP-norm} implies \eqref{LP-oper}. Other properties are established similarly  using \eqref{LP-split}--\eqref{LP-norm}.
      \end{proof}
%\begin{lem}
\end{appendices}
\section{Acknowledgments}
This research did not receive any specific grant from funding agencies in the public, commercial, or
not-for-profit sectors.
%\end{lem}
    %   $\sum_{j=0}^\infty P_j E_B(\lambda) = E_\lambda$ (we kill off high eigenvalues/frequencies of our operator).\\
      % For simplicity, given $\triangle R>0$, we let $R_1=R-\triangle R$. Define $u_R = E_B(R-\triangle R,R)u$ Note $A u_R =\int_{R-\triangle R }^R\lambda dE_B(\lambda)u $.
    %   \hide{Essentially $u_R$ is $d E_R(u)$ - it is convenient to do projections twice}
%\section{Acknowledgments}
%The authors would like to express gratitude to Cristian Rios, who introduced us to the problem of hypoellipticity and initiated the program that lead to this result.
\bibliographystyle{amsalpha}
\bibliography{Hypo-paper-2022}
\end{document}